\documentclass[10pt]{article}
\usepackage[english]{babel}
\usepackage{geometry}
\usepackage{amssymb}
\usepackage{indentfirst}
\usepackage{amsmath}
\usepackage{amsthm}
\newtheorem{proposition}{Proposition}[section]
\newtheorem{lemma}{Lemma}[section]
\newtheorem{theorem}{Theorem}[section]
\usepackage{mathrsfs}

\title{Proof of Quasi Riemann Hypothesis using $\mathbb Q(i)$}
\author{XinHang Ji}
\date{}

\begin{document}
	\maketitle
	
	\begin{abstract}
		Using quartic Gauss sums in place of the cubic Gauss sums in
		OpenAI's original proof, we prove that the Dedekind $\zeta$-function
		$\zeta_{\mathbb Q(i)}(s)$ of $\mathbb Q(i)$ has no zeros in the half-plane
		$\Re s>599/600$, thereby obtaining the same zero-free half-plane for
		the Riemann $\zeta$-function.
		We obtain the upper bound for the smoothed sum by applying Poisson
		summation twice and using the functional equation of the quartic
		Kubota series of angular weight $1$.
		In the Mellin representation, a zero-density bound for the family of
		quartic characters controls the terms from nonprincipal characters,
		and a residue and a lattice average over fourth powers give the
		main term.
	\end{abstract}
	
	\section{Introduction}

	This paper aims to provide a shorter, more readable proof of the
	Quasi-Riemann Hypothesis.
	
	The Quasi-Riemann Hypothesis asserts that the Riemann $\zeta$-function
	has no zeros in a right half-plane whose boundary has real part
	strictly less than $1$.
	The two OpenAI papers establish zero-free half-planes with boundaries
	at $7/8$ and $11/12$, respectively \cite{qrh,qrh-alt}, using structure
	unique to the number field $\mathbb Q(\sqrt{-3})$.
	Here we take a different route that dispenses with $\mathbb Q(\sqrt{-3})$
	and uses quartic Gauss sums over the Gaussian integer ring $\mathbb Z[i]$
	to study the Dedekind $\zeta$-function
	$\zeta_{\mathbb Q(i)}(s)=\zeta(s)L(s,\chi_{-4})$ and prove the following
	theorem.
	
	\begin{theorem}\label{thm:main}
		$\zeta_{\mathbb Q(i)}(s)$ has no zeros in $\Re s>599/600$.
	\end{theorem}
	
	Since $L(s,\chi_{-4})$ is entire, the theorem also shows that the
	Riemann $\zeta$-function has no zeros in the same half-plane.
	
	Three features account for most of the reduction in length:
	
	1. The quartic theta functions used by David, Dunn, Hamieh, and Lin
	come from residues of Eisenstein series of angular weight $0$
	\cite[Section~4.5]{ddhl}.
	In estimating the smoothed sum, this paper uses the functional equation
	of the quartic Kubota series of angular weight $1$.
	At this angular weight, the completed series at the cusps are entire,
	so the contour shift for the reflection picks up no residues
	\cite[Proposition~4.3]{ddhl}.
	The angular factor after reflection makes the zero-frequency term
	in the second application of Poisson summation vanish.
	
	2. For the nonprincipal contribution, we group the characters by zero
	location and use de Faveri's zero-density bound for the family of
	quartic characters directly to control the number of characters in each group
	\cite[Corollary~1.5]{defaveri}, dispensing with separate derivations
	of zero-detection and large-sieve estimates.
	
	3. This paper limits its treatment of the family to one fixed
	Dedekind $\zeta$-function.
	
	\section{Notation and main propositions}\label{sec:statements}
	
	\subsection{Arithmetic notation}
	Throughout this paper, elements of $\mathbb Z$ are called integers
	and elements of $\mathbb Z[i]$ are called Gaussian integers.
	A prime is a positive prime in the integers, and a Gaussian prime is
	a prime element of the Gaussian integer ring.
	Let $\mu_K$ denote the M\"obius function on the Gaussian integer ring
	and $\mu$ the usual M\"obius function on the positive integers.
	
	Write $\lambda=1+i$.
	The norm of a nonzero integral ideal $\mathfrak a$ of $\mathbb Z[i]$
	is denoted by $N\mathfrak a=\#(\mathbb Z[i]/\mathfrak a)$.
	If $\mathfrak a=(a)$ is a nonzero integral ideal, then
	$N\mathfrak a=|a|^2$.
	A Gaussian integer or ideal coprime to $\lambda$ is called odd.
	The ring $\mathbb Z[i]$ is a Euclidean domain, and its four units
	represent all invertible residue classes modulo $\lambda^3$, so every
	nonzero odd ideal has a unique generator $a\equiv1\pmod{\lambda^3}$,
	called its primary generator.
	Sums over odd ideals are always taken over their primary generators,
	and the ranges of sums over lattice points or residue classes are
	specified explicitly.
	For a nonzero complex number $a$, set $\alpha(a)=a/|a|$.
	For a Gaussian prime $\tau$, let $v_\tau$ denote the corresponding
	valuation on $\mathbb Q(i)$, with $v_\tau(\tau)=1$ and
	$v_\tau(0)=+\infty$.
	For a nonzero Gaussian integer $a$ and a nonnegative integer $k$,
	write $\tau^k\Vert a$ for $v_\tau(a)=k$.
	Let $\omega(a)$ count the distinct Gaussian prime factors of $a$ up
	to associates.
	The function $\mu_K$ takes the value $(-1)^{\omega(a)}$ on nonzero
	squarefree Gaussian integers $a$ and zero on all other Gaussian integers.
	The indicator of the condition $A$ is denoted by $\mathbf 1_A$.
	For positive $x,U$, write $x\sim U$ to mean $U\le x<2U$, and write
	$x\asymp U$ to mean that their ratio is bounded above and below by
	positive constants.
	The notation $A\ll B$ means $|A|\le C B$, where the constant $C$ is
	independent of the summation scales.
	The constants in $\ll$ and $O$ may depend on the fixed smooth weights
	and on the parameters listed in subscripts.
	
	$L(s,\xi)$ denotes the $L$-function of the Hecke character $\xi$, and
	$L_\lambda(s,\xi)$ denotes this function with the factor
	$(1-\xi(\lambda)|\lambda|^{-2s})^{-1}$ omitted.
	We have $\zeta_{\mathbb Q(i)}(s)=\sum_{\mathfrak a\ne0}(N\mathfrak a)^{-s}$
	($\Re s>1$), and in particular
	$\zeta_\lambda(s)=L_\lambda(s,1)=(1-2^{-s})\zeta_{\mathbb Q(i)}(s)$.
	
	For an integer $m$, $\alpha^{4m}$ denotes the angular character on
	ideals given by $(a)\mapsto \alpha(a)^{4m}$.
	Its value is independent of the choice of unit multiplying the generator.
	Let $\tau$ be an odd Gaussian prime and $x$ a Gaussian integer.
	If $\tau\nmid x$, define $\left(\frac{x}{\tau}\right)_4$ as the fourth root of unity
	congruent to $x^{(|\tau|^2-1)/4}$ modulo $\tau$.
	The symbol is zero when $\tau\mid x$.
	For odd Gaussian integers $c$, extend the symbol $\left(\frac{x}{c}\right)_4$
	multiplicatively over the Gaussian prime factors of $c$, with the
	convention $\left(\frac{x}{1}\right)_4=1$.
	For primary $u$, let $\xi_u$ denote the primitive Hecke character
	underlying the character $(a)\mapsto\left(\frac{u}{a}\right)_4$ on ideals coprime to
	$\lambda u$.
	In the unique factorization $u=u_0f^4$, the factors $u_0,f$ are primary
	and $0\le v_\tau(u_0)\le3$.
	On ideals coprime to $\lambda u$, the fourth-power factor takes the
	value $1$, so $\xi_u=\xi_{u_0}$.
	By \cite[Example~2.3 and Section~2.5.1]{defaveri}, its conductor norm
	is $O(|u_0|^2)$, and every Gaussian prime factor of $u_0$ occurs in
	the primitive conductor.
	Thus $\xi_u$ is principal if and only if $u_0=1$.
	
	\subsection{Decomposition of the smoothed sum}
	
	For a complex number $x$, put $e_0(x)=\exp(2\pi i\Re x)$.
	Let $r$ be a Gaussian integer and let the modulus $c$ be an odd Gaussian
	integer.
	The following formulas define two Gauss sums, where the Gaussian integer
	$x$ ranges over a complete residue system modulo $c$.
	\[
	g_+(r,c)=\sum_{x\bmod c}\left(\frac{x}{c}\right)_4e_0(rx/c),
	\qquad
	g_-(r,c)=\sum_{x\bmod c}\overline{\left(\frac{x}{c}\right)_4}e_0(rx/c).
	\]
	For $\nu\in\lambda^{-2}\mathbb Z[i]$, write $g_4(\nu,c)=g_+(2\nu,c)$
	for the quartic Gauss sum with the trace additive character.
	For primary $c,b$, put $U(c,b)=(-1)^{((|c|^2-1)/4)((|b|^2-1)/4)}$.
	The following equation defines a periodic function $F_{c,b}(u)$ on the
	Gaussian integers by its finite Fourier transform modulo $cb$.
	\begin{equation}\label{eq:fourier}
		\sum_{u\bmod cb}F_{c,b}(u)e_0(ru/(cb))
		=\alpha(c)\left(\frac{-1}{c}\right)_4U(c,b)
		\overline{g_+(r^3,c)}\,g_-(r,b).
	\end{equation}
	Equation~\eqref{eq:fourier} holds for every $r\bmod cb$.
	The pairing $(r,u)\mapsto e_0(ru/(cb))$ is nondegenerate, so finite
	Fourier inversion determines $F_{c,b}$ uniquely.
	
	The function $J$ modulo $\lambda^4$ takes the values $1$ and $-1$
	on the classes of $1$ and $1+\lambda^3$, respectively, and is zero
	on all other classes.
	These two classes form the group of primary residue classes, so $J$
	is multiplicative on this group and $J(f^4)=1$ for primary $f$.
	
	Let $C_c^\infty((1,2))$ denote the space of smooth functions with
	compact support in $(1,2)$.
	In this space, fix nonnegative and nonzero $V,W_1$ and take an arbitrary,
	possibly complex-valued $W$.
	The weights are zero outside their supports.
	Take real parameters $X\ge1$, $1\le H\le X$, and $Y\ge1$, and a
	primary $d$. Define
	\begin{equation}\label{eq:Bd}
		\mathcal B_d(u)=\sum_{c,b} F_{c,b}(u)
		W\left(\frac{|c|^2\,|d|^8}{X}\right)W_1\left(\frac{|b|^2}{Y}\right).
	\end{equation}
	The smoothed sum under consideration is
	\begin{equation}\label{eq:I}
		I(X,H,Y)=\sum_d|d|^6\alpha(d)^4
		\sum_{u\in\mathbb Z[i]}J(u)V\left(\frac{|u|^2\,|d|^8}{H}\right)\mathcal B_d(u).
	\end{equation}
	Since $J$ is supported on primary Gaussian integers, the factorization
	$u=u_0 f^4$ gives $I=I_{\mathrm{np}}+I_{\mathrm{pr}}$, where
	\begin{equation}\label{eq:I-parts}
		\begin{aligned}
			I_{\mathrm{np}}(X,H,Y)
			&=\sum_d |d|^6\alpha(d)^4
			\sum_{\substack{u\in 1+\lambda^3\mathbb Z[i]\\ u_0\ne1}}
			J(u)V\left(\frac{|u|^2\,|d|^8}{H}\right)\mathcal B_d(u),\\
			I_{\mathrm{pr}}(X,H,Y)
			&=\sum_d |d|^6\alpha(d)^4
			\sum_f V\left(\frac{|f|^8\,|d|^8}{H}\right)\mathcal B_d(f^4).
		\end{aligned}
	\end{equation}
	In the formulas above, the summation indices $c,b,d,f$ all run over
	primary Gaussian integers.
	Here $\mathrm{np}$ and $\mathrm{pr}$ stand for nonprincipal and principal
	characters, respectively.
	The support of $V$ restricts both sums to $|d|^2\ll H^{1/4}$.
	
	\subsection{Mellin transform of the main term}\label{sec:mellin-notation}
	
	We define the main term for the principal character in the half-plane
	$\Re z>1$.
	For an odd Gaussian prime $\tau$, write $q=|\tau|^2$.
	For a nonnegative integer $j$, define the following local factors,
	where $k,\ell$ range over the nonnegative integers.
	\begin{equation}\label{eq:local-D}
		D_{\tau,j}(z)=\sum_{k,\ell\ge0}
		\frac{F_{\tau^k,\tau^\ell}(\tau^{4j})}{q^{kz+\ell}},
		\qquad
		P_{\tau,j}(z)=\frac{1-q^{-1}}{1-q^{-z}}D_{\tau,j}(z).
	\end{equation}
	Averaging the local factors using the lattice densities for
	$v_\tau(f)=j$ gives
	\begin{equation}\label{eq:mean-euler}
		\mathcal Q(z)=\prod_{\tau\nmid\lambda}
		\left((1-q^{-1})\sum_{j\ge0}q^{-j}P_{\tau,j}(z)\right).
	\end{equation}
	Lemma~\ref{lem:nonzero} will show that the series and product above
	converge absolutely in $\Re z>1-2^{-6}$ and that $\mathcal Q$ is
	holomorphic and nonzero in that half-plane.
	
	For a smooth, compactly supported function $h$ on $(0,\infty)$, denote
	its Mellin transform by $\widehat h(z)=\int_0^\infty h(x)x^{z-1}\,\mathrm dx$,
	and for real $\sigma$ write
	$\int_{(\sigma)}=\int_{\sigma-i\infty}^{\sigma+i\infty}$.
	The function $S_W(X)$ appearing in the main term is the inverse Mellin
	transform of $\widehat W(s)D(s)$, where $D$ and $S_W$ are defined by
	\begin{equation}\label{eq:mellin-S}
		D(z)=\frac{L_\lambda(4z-2,\alpha^4)\mathcal Q(z)}{\zeta_\lambda(z)},
		\qquad
		S_W(X)=\frac{1}{2\pi i}\int_{(2)}
		\widehat W(s)X^sD(s)\,\mathrm ds \quad (X>0).
	\end{equation}
	Absolute convergence of the Euler product and decay of $\widehat W$
	allow the integral to be taken on any line with $\Re s>1$, so
	$S_W(X)\ll_{\varepsilon,W} X^{1+\varepsilon}$ for every $\varepsilon>0$.
	On $\Re s=a\ge2$, integration by parts gives
	$S_W(X)\ll_W(2X)^a(1+a)^2$.
	Letting $a\to\infty$ gives $S_W(X)=0$ ($0<X<1/2$).
	Mellin inversion gives
	\begin{equation}\label{eq:mellin-identity}
		\widehat W(z)D(z)=\int_0^\infty S_W(X)X^{-z}\,\frac{\mathrm dX}{X}
		\qquad(\Re z>1).
	\end{equation}
	
	The main-term constant is
	\[
	K_0=\left(\frac{\pi}{8}\right)^2
	\widehat W_1(1)\int_0^\infty V(x^4)\,\mathrm dx>0.
	\]
	
	\subsection{Three propositions}
	
	For the estimates of the principal and nonprincipal contributions, fix
	\begin{equation}\label{eq:parameters}
		H=X^{65/144},\qquad Y=X^{1/288} .
	\end{equation}
	
	\begin{proposition}[Upper bound for the smoothed sum]\label{prop:dual}
		For $X\ge1$, $1\le H\le X$, $Y\ge1$, and any $\varepsilon>0$,
		\begin{equation}\label{eq:dual-bound}
			|I(X,H,Y)|\ll_\varepsilon (XHY)^\varepsilon X^2Y^5H^{-2}.
		\end{equation}
	\end{proposition}
	
	\begin{proposition}[Estimate of the nonprincipal contribution]\label{prop:nonprincipal}
		With the parameters in \eqref{eq:parameters}, for any $\varepsilon>0$,
		\begin{equation}\label{eq:nonprincipal}
			I_{\mathrm{np}}(X,H,Y)
			\ll_\varepsilon X^{599/600+\varepsilon}YH^{1/4}.
		\end{equation}
	\end{proposition}
	
	\begin{proposition}[Estimate of the main term]\label{prop:principal}
		With the parameters in \eqref{eq:parameters}, for any $\varepsilon>0$,
		\begin{equation}\label{eq:principal}
			\begin{aligned}
				I_{\mathrm{pr}}(X,H,Y)
				&=K_0YH^{1/4}S_W(X)\\
				&\quad+O_\varepsilon\left(X^{1+\varepsilon}Y^{13/25}H^{1/4}
				+X^{1+\varepsilon}YH^{1/8}\right).
			\end{aligned}
		\end{equation}
	\end{proposition}
	
	\subsection{Outline of the proof}\label{sec:proof-outline}
	
	Proposition~\ref{prop:dual} directly estimates the whole smoothed sum $I$.
	After applying Poisson summation to the sum over $u$, denote the Fourier
	frequency by $r$.
	Insert the resulting Gauss sums into the quartic Kubota series of angular
	weight $1$ and apply its functional equation to obtain the angular factor
	$\alpha(r)^3$ after reflection.
	The angular factor has zero average over the circle, so the zero-frequency
	term vanishes when we apply Poisson summation again in $r$.
	Bounds for finite Gauss sums and the decay of the smooth weights control
	the remaining frequencies.
	
	Propositions~\ref{prop:nonprincipal} and \ref{prop:principal} use the
	decomposition \eqref{eq:I-parts} to treat the same $I$.
	In $\mathcal B_d(u)$, apply Mellin inversion separately to the norms of
	$c,b$, take $z,w$ as the corresponding integration variables, and
	factor out
	$L_\lambda(w,\overline{\xi_u})/L_\lambda(z,\xi_u)$.
	For the nonprincipal contribution, group the characters by the largest
	real part among their zeros up to the truncation height.
	For the group in which this real part is at least $599/600$, keep the $z$
	contour fixed and move the $w$ contour to $\Re w=13/25$.
	For the remaining groups, keep the $w$ contour fixed, move the $z$
	contour, and use the zero-density bound to control the contribution
	of each group.
	For the principal contribution, the residue at $w=1$, averaged over
	the lattice points $f$, gives $K_0YH^{1/4}S_W(X)$, with the local
	averages combining to form $\mathcal Q(z)$.
	
	Combining the three propositions gives
	$S_W(X)\ll_{\varepsilon,W} X^{599/600+\varepsilon}$.
	The integral in \eqref{eq:mellin-identity} therefore extends
	$\widehat W(z)D(z)$ holomorphically to $\Re z>599/600$.
	Choosing at each point a weight with $\widehat W\ne0$ and then dividing
	by the nonzero numerator of $D$ gives Theorem~\ref{thm:main}.
	
	\section{Proof of the upper bound for the smoothed sum}\label{sec:dual}
	
	\subsection{Local estimates for quartic coefficients}\label{sec:local}
	
	\begin{lemma}[Square of the Gauss sum]\label{lem:gauss-square}
		If $c$ is squarefree, then
		\[
		\frac{g_+(1,c)^2}{|c|^2}=\mu_K(c)\alpha(c)\left(\frac{-1}{c}\right)_4,
		\qquad F_{c,1}(u)=\mu_K(c)\left(\frac{u}{c}\right)_4.
		\]
		For any primary $b$, we have $F_{1,b}(u)=\overline{\left(\frac{u}{b}\right)_4}$.
	\end{lemma}
	
	\begin{proof}
		First let $\tau=a+bi$ be a split Gaussian prime, with $a,b$ integers
		and $p=|\tau|^2$ prime.
		Then $a$ is odd and $\bar\tau\equiv 2a\pmod\tau$.
		Quadratic reciprocity gives $(a/p)=(p/|a|)=1$, since
		$p\equiv b^2\pmod{|a|}$.
		Therefore
		\[
		\left(\frac{\bar\tau}{\tau}\right)_4^2
		=\left(\frac{2}{p}\right)=\left(\frac{-1}{\tau}\right)_4.
		\]
		The Gauss sum $g_4$, defined using the trace additive character,
		satisfies $g_4(1,\tau)^2=-\tau\sqrt p$ (see \cite[Eq.~(3.11)]{ddhl}).
		The desired formula for the square follows from
		$g_+(1,\tau)=g_4(1/2,\tau)=\left(\frac{2}{\tau}\right)_4g_4(1,\tau)$ and
		$\left(\frac{2}{\tau}\right)_4^2=\left(\frac{-1}{\tau}\right)_4$.
		
		If $\tau=-p$, where the prime $p$ is inert in $\mathbb Q(i)$,
		then $\left(\frac{\cdot}{\tau}\right)_4$ is trivial on $\mathbb F_p^\times$.
		Sum over the fibers of the trace map $\mathbb F_{p^2}\to\mathbb F_p$.
		The fiber over zero contributes $(p-1)\left(\frac{i}{\tau}\right)_4$, and the fibers
		over nonzero values all contribute equally, so the vanishing of the
		total sum gives a contribution of $-\left(\frac{i}{\tau}\right)_4$ from each fiber
		over a nonzero value.
		Thus $g_+(1,\tau)=p\left(\frac{i}{\tau}\right)_4$.
		From $p^2\equiv1\pmod8$ we have $\left(\frac{i}{\tau}\right)_4^2=1$.
		This gives $g_+(1,\tau)^2/|\tau|^2=1=\mu_K(\tau)\alpha(\tau)\left(\frac{-1}{\tau}\right)_4$.
		
		For coprime primary $c_1,c_2$, the Chinese remainder theorem gives
		\[
		g_+(r,c_1c_2)
		=\left(\frac{c_2}{c_1}\right)_4\left(\frac{c_1}{c_2}\right)_4g_+(r,c_1)g_+(r,c_2).
		\]
		By quartic reciprocity \cite[Eq.~(3.2)]{ddhl}, the square of the cross
		factor is $1$.
		Multiplying the identities for the prime factors gives the formula
		for squarefree $c$.
		In this case $|g_+(1,c)|=|c|$, and when $(r,c)=1$ we have
		$g_+(r^3,c)=\left(\frac{r}{c}\right)_4g_+(1,c)$, so
		\[
		\alpha(c)\left(\frac{-1}{c}\right)_4\overline{g_+(r^3,c)}
		=\mu_K(c)g_+(r,c).
		\]
		If $(r,c)\ne1$, both sides are zero.
		Finite Fourier inversion of \eqref{eq:fourier} gives the formulas
		for $F_{c,1}$ and $F_{1,b}$.
	\end{proof}
	
	\begin{lemma}[Scaling and multiplicativity]\label{lem:crt}
		Let $t,u$ be Gaussian integers.
		If $(t,cb)=1$, then
		$F_{c,b}(tu)=\left(\frac{t}{c}\right)_4\overline{\left(\frac{t}{b}\right)_4}F_{c,b}(u)$.
		If $(c_1b_1,c_2b_2)=1$, then
		\begin{equation}\label{eq:crt}
			F_{c_1c_2,b_1b_2}(u)=F_{c_1,b_1}(u)F_{c_2,b_2}(u).
		\end{equation}
		Moreover, $F_{c,b}(u)$ vanishes unless $(c,b)\mid u$.
	\end{lemma}
	
	\begin{proof}
		Replacing $r$ by $tr$ on the right-hand side of \eqref{eq:fourier}
		allows us to factor out $\overline{\left(\frac{t}{c}\right)_4}\left(\frac{t}{b}\right)_4$.
		The inverse transform gives the stated scaling formula.
		
		By the Chinese remainder theorem, the two arguments for $u$ in the
		inverse transform become $(c_2b_2)^{-1}u\bmod c_1b_1$ and
		$(c_1b_1)^{-1}u\bmod c_2b_2$, respectively.
		Of the cross factors from the Gauss sum multiplication formula and the
		scaling formula, the part involving only $c$ is
		$\overline{\left(\frac{c_2}{c_1}\right)_4\left(\frac{c_1}{c_2}\right)_4}^{\,2}=1$, and the part
		involving only $b$ cancels.
		The mixed part is
		\[
		\overline{\left(\frac{b_2}{c_1}\right)_4}\left(\frac{c_1}{b_2}\right)_4
		\,\overline{\left(\frac{b_1}{c_2}\right)_4}\left(\frac{c_2}{b_1}\right)_4
		=U(c_1,b_2)U(c_2,b_1).
		\]
		The function $U$ is completely multiplicative in each variable, so this
		factor cancels the two mixed factors in the expansion of $U(c_1c_2,b_1b_2)$,
		giving \eqref{eq:crt}.
		Finally, the period in $r$ of the right-hand side of \eqref{eq:fourier}
		divides $\operatorname{lcm}(c,b)$, so its inverse transform vanishes
		unless $(c,b)\mid u$.
	\end{proof}
	
	\begin{lemma}[Coefficients at Gaussian prime powers]\label{lem:prime-power}
		Let $\tau$ be an odd Gaussian prime and $u$ a nonzero Gaussian integer,
		and write $q=|\tau|^2$ and $v=v_\tau(u)$.
		For any nonnegative integers $k,\ell$,
		\begin{equation}\label{eq:local}
			|F_{\tau^k,\tau^\ell}(u)|
			\le 4q^{\min(v,k,\lfloor2k/3\rfloor+\ell)}.
		\end{equation}
		When $k=2,3$, the coefficients are identically zero.
		If $v=0$, the only possibly nonzero terms are those with $k=0$ or
		$(k,\ell)=(1,0)$.
		If $v>0$, the terms with $k=0,\ell>0$ and $(k,\ell)=(1,0)$ are zero,
		and for $k=1,\ell>0$ only the term with $\ell=1$ can be nonzero,
		with $|F_{\tau,\tau}(u)|\le q-1$.
	\end{lemma}
	
	\begin{proof}
		For $k>0$ with $4\nmid k$, the conductor of $\left(\frac{\cdot}{\tau^k}\right)_4$ is $\tau$.
		Summing over the fibers of reduction modulo $\tau$ shows that for a
		Gaussian integer $r$, the sum $g_+(r^3,\tau^k)$ can be nonzero only
		if there are a nonnegative integer $j$ and a Gaussian integer $x$
		such that $k=3j+1$, $r=\tau^j x$, and $\tau\nmid x$, in which case it
		is the product of a nontrivial character modulo $\tau$ and a constant
		of modulus $q^{k-1/2}$.
		For $k>0$ with $4\mid k$, it equals
		\[
		q^k\mathbf 1_{\tau^{\lceil k/3\rceil}\mid r}
		-q^{k-1}\mathbf 1_{\tau^{\lceil(k-1)/3\rceil}\mid r}.
		\]
		The same formulas hold for $g_-(r,\tau^\ell)$, with the nontrivial
		character appearing when $\ell=j+1$ and the two exponents in the
		divisibility conditions for the trivial character given by $\ell$ and
		$\ell-1$.
		The Gauss sum with exponent zero is the constant $1$.
		
		Substitute these two types of expressions into \eqref{eq:fourier}.
		If the supports of the two nontrivial terms intersect, then $k=3j+1$
		and $\ell=j+1$, with $4\nmid k$ and $4\nmid\ell$.
		If $j\equiv1\pmod4$, then $4\mid k$, and if $j\equiv3\pmod4$, then
		$4\mid\ell$.
		Thus $j$ is even and $k\equiv\ell\pmod4$.
		With $r=\tau^j x$, the characters in the conjugated first Gauss sum
		and the second Gauss sum are $\left(\frac{x}{\tau}\right)_4^{-k}$ and $\left(\frac{x}{\tau}\right)_4^\ell$,
		respectively, and their product is $1$ when $\tau\nmid x$.
		In this case, write the indicator of the unit condition as
		$\mathbf 1_{\tau^j\mid r}-\mathbf 1_{\tau^{j+1}\mid r}$.
		The Gauss sum product thus has at most four components.
		The inverse transform modulo $\tau^{k+\ell}$ takes
		$A\mathbf 1_{\tau^j\mid r}$ to a function supported on
		$v\ge k+\ell-j$ with modulus at most $|A|q^{-j}$.
		If the component also carries a nontrivial character $\psi(r/\tau^j)$,
		its inverse transform is supported on $v=k+\ell-j-1$ and has modulus
		at most $|A|q^{-j-1/2}$.
		These two conclusions follow from orthogonality of additive characters
		and the Gauss sum formula modulo $\tau$, respectively.
		
		For a constant component, $|A|\le q^{k+\ell}$, and when the two
		nontrivial characters cancel, we also have $|A|\le q^{k+\ell-1}$.
		For a component with a nontrivial character, $|A|\le q^{k+\ell-1/2}$.
		Thus the support conditions above give the bound $q^v$.
		If the component of $g_-$ is constant, $q^{-j}$ absorbs its modulus.
		If the component carries a nontrivial character, then $j=\ell-1$.
		The factor $q^{-1/2}$ in the inverse transform or the extra $q^{-1/2}$
		in the first Gauss sum completes the normalization.
		This gives the bound $q^k$.
		For the two constant components of the first Gauss sum, use
		$k-\lceil k/3\rceil=\lfloor2k/3\rfloor$ and
		$k-1-\lceil(k-1)/3\rceil\le\lfloor2k/3\rfloor$, respectively.
		If the first character is nontrivial, then $k=3j+1$, so
		$k-1-j=\lfloor2k/3\rfloor$, with the extra half power absorbed by the
		inverse transform or the second nontrivial Gauss sum.
		This gives the third bound $q^{\lfloor2k/3\rfloor+\ell}$, and summing
		the components proves \eqref{eq:local}.
		
		For $k=2,3$, the equation $k=3j+1$ has no nonnegative integer solution.
		If $v=0$ and $k,\ell>0$, the support condition in Lemma~\ref{lem:crt}
		forces the coefficient to vanish.
		For $k>1,\ell=0$, the support condition for the inverse transform above forces
		$v>0$.
		If $v>0$, the claims for $k=0$ or $(k,\ell)=(1,0)$ follow from
		Lemma~\ref{lem:gauss-square}.
		For $k=1,\ell>0$, the first Gauss sum can be nonzero only if
		$\tau\nmid r$, so only $\ell=1$ can give a nonzero term.
		In this case the two characters cancel, and the inverse transform is
		a Ramanujan sum modulo $\tau^2$ times a constant of modulus $q^{-1}$,
		so the inverse transform has modulus at most $q-1$.
	\end{proof}
	
	\subsection{First Poisson summation and quartic reflection}
	
	The Fourier frequency $r$ is a Gaussian integer.
	The finite Fourier transform of $J$ is
	\begin{equation}\label{eq:j-transform}
		\widehat J(r)=\sum_{v\bmod\lambda^4}J(v)e_0(rv/\lambda^4)
		=e_0(r/\lambda^4)\bigl(1-e_0(r/\lambda)\bigr).
	\end{equation}
	Since $e_0(r/\lambda)=1$ if and only if $\lambda\mid r$, the transform
	$\widehat J$ is supported only on odd $r$.
	The planar Fourier transform uses the kernel $e_0(-x\xi)$ and planar
	area measure $\mathrm dA(x)$.
	For complex $\xi$, define the radial function $\mathcal F_V$ by
	\[
	\mathcal F_V(|\xi|^2)=\int_{\mathbb C}V(|x|^2)e_0(-x\xi)\,\mathrm dA(x).
	\]
	The function and all its radial derivatives decay rapidly.
	
	\begin{lemma}[First Poisson summation]\label{lem:poisson}
		For primary $b$, put $A_b=X|b|^2/H$.
		For real $v\ge0$ and $x>0$, define
		\[
		w_{v,b}(x)=W(x)x^{-1/2}\mathcal F_V\left(\frac{v}{|\lambda|^8\,A_bx}\right).
		\]
		For $\beta\in\{1,1+\lambda^3\}$, put
		\begin{equation}\label{eq:column}
			\mathcal S_\beta(r;b)
			=\sum_{\substack{c\equiv\beta\bmod\lambda^4\\ d}}
			|d|^2\,\alpha(d)^4\alpha(c)\frac{\overline{g_+(r^3,c)}}{|c|}\,
			w_{|r|^2,b}\left(\frac{|c|^2\,|d|^8}{X}\right).
		\end{equation}
		With $C_\beta(b)=J(\beta b)\left(\frac{-1}{\beta}\right)_4U(\beta,b)$, we have
		\begin{equation}\label{eq:poisson}
			I=\frac{H}{|\lambda|^8\sqrt X}\sum_b W_1\left(\frac{|b|^2}{Y}\right)
			\sum_\beta C_\beta(b)
			\sum_{\lambda\nmid r}\widehat J(r)\frac{g_-(r,b)}{|b|^2}
			\mathcal S_\beta(r;b).
		\end{equation}
	\end{lemma}
	
	\begin{proof}
		Apply Poisson summation to the sum over $u$ in \eqref{eq:I}, using
		residue classes modulo $\lambda^4cb$.
		Writing $u=(cb)u_1+\lambda^4u_2$ gives $J((cb)u_1)=J(cb)J(u_1)$.
		By the scaling formula in Lemma~\ref{lem:crt},
		$F_{c,b}(\lambda^4u_2)=F_{c,b}(u_2)$.
		Thus the resulting finite Fourier sum is the right-hand side of
		\eqref{eq:fourier} multiplied by $J(cb)\widehat J(r)$.
		The argument and normalization factor of the continuous transform
		are, respectively,
		\[
		\frac{|r|^2X}{|\lambda|^8\,A_b\,|c|^2\,|d|^8},
		\qquad \frac{H}{|d|^8\,|\lambda|^8\,|c|^2\,|b|^2}.
		\]
		Multiplying the normalization factor by $|d|^6$ and factoring
		$|c|$ out of the Gauss sum gives
		\[
		\frac{H}{|\lambda|^8\,|b|^2\,|c|\,|d|^2}.
		\]
		We have $\left(\frac{-1}{c}\right)_4=(-1)^{(|c|^2-1)/4}$, and $|c|^2\bmod 8$ is determined
		by $c\bmod\lambda^4$.
		Therefore, when $c\equiv\beta\bmod\lambda^4$,
		\[
		J(cb)\left(\frac{-1}{c}\right)_4U(c,b)=C_\beta(b),\qquad |C_\beta(b)|\le1.
		\]
		The Poisson summation formula applies for all $d$, and the identity
		still holds when the original sum is empty.
		The support of $W$ makes the ranges of $c,d$ finite, and the smooth
		weight ensures absolute convergence of the sum over $r$.
	\end{proof}
	
	\begin{lemma}[Quartic reflection formula]\label{lem:reflection}
		Fix $\beta\in\{1,1+\lambda^3\}$.
		In the formula below, the integer $j$ indexes finitely many terms
		in the cusp expansions.
		In each term, $\beta_j\in\{1,1+\lambda^3\}$, $2\kappa_j$ is a unit
		times a nonnegative integer power of $\lambda$, and $\eta_j(r)$ is
		a bounded function periodic modulo $\lambda^7$.
		For $\lambda\nmid r$,
		\begin{equation}\label{eq:reflection}
			\begin{aligned}
				\mathcal S_\beta(r;b)
				={}&\frac{X\,\alpha(r)^3}{|r|^3}
				\sum_j \eta_j(r)
				\sum_{\substack{n\equiv\beta_j\bmod\lambda^4\\ e}}
				|e|^2\,\overline{\alpha(ne^4)}
				\frac{\overline{g_4(\kappa_j r^3,n)}}{|n|}\\
				&\times
				\mathcal W_{j,b}\left(|r|^2,\frac{|n|^2\,|e|^8X}{|r|^6}\right).
			\end{aligned}
		\end{equation}
		Here $n,e$ are primary.
		The kernel $\mathcal W_{j,b}(v,y)$ is defined below.
		For nonnegative integers $a_1,a_2$ and real $B,L>0$,
		\begin{equation}\label{eq:kernel}
			\left|(v\partial_v)^{a_1}(y\partial_y)^{a_2}
			\mathcal W_{j,b}(v,y)\right|
			\ll (1+v/A_b)^{-L}
			\begin{cases}
				y^{-9/8},&0<y\le1,\\
				y^{-B},&y\ge1.
			\end{cases}
		\end{equation}
	\end{lemma}
	
	\begin{proof}
		Putting $\nu=r^3/2$ gives $g_4(\nu,c)=g_+(r^3,c)$ and $v_\lambda(\nu)=-2$.
		Use the normalization in \cite[Eqs.~(4.1) and (4.7)]{ddhl} to define
		the series below for angular weight $\ell\in\mathbb Z$.
		\[
		\Psi_\beta(s,\nu,\ell)=
		\sum_{c\equiv\beta\bmod\lambda^4}
		\frac{g_4(\nu,c)\overline{\alpha(c)}^{\,\ell}}{|c|^{2s}},
		\qquad
		Z_\beta(s,\nu,\ell)=L_\lambda(4s-3,\alpha^{-4\ell})
		\Psi_\beta(s,\nu,\ell).
		\]
		These series converge absolutely in $\Re s>3/2$.
		Apply Mellin inversion to \eqref{eq:column}.
		Since the sum over $d$ gives $L_\lambda(4z-1,\alpha^4)$, we obtain
		\[
		\mathcal S_\beta(r;b)=\frac{1}{2\pi i}\int_{(9/8)}
		\widehat w_{|r|^2,b}(z)X^z\,\overline{Z_\beta(\bar z+1/2,\nu,1)}\,\mathrm dz.
		\]
		
		Use the cusp numbering $i=1,\ldots,24$ of
		\cite[Eqs.~(4.11) and (4.24)--(4.28)]{ddhl}.
		Let $\psi^{(4)}_{i1}(s,\nu,\ell)$ be the cusp series in Eq.~(4.11)
		of that paper, and put
		\[
		Z_{i1}(s,\nu,\ell)=L_\lambda(4s-3,\alpha^{-4\ell})
		\psi^{(4)}_{i1}(s,\nu,\ell),
		\qquad \widehat Z_{i1}(s,\nu,\ell)=G_\infty(s,\ell)Z_{i1}(s,\nu,\ell),
		\]
		The archimedean factor is
		\[
		G_\infty(s,\ell)=8(2\pi)^{-3s-\frac32|\ell|+\frac32}
		\prod_{m=1}^3\Gamma\left(s+\frac{|\ell|}{2}-\frac m4\right).
		\]
		Here $\Gamma$ is the gamma function.
		Let $i_\beta$ be the index of cusp $0$ when $\beta=1$ and of cusp
		$\lambda^4/(1+\lambda^3)$ when $\beta=1+\lambda^3$.
		By Eqs.~(4.12)--(4.13) of that paper,
		$Z_\beta(s,\nu,1)=-J(\beta)\,|\lambda|^8\,Z_{i_\beta1}(s,\nu,1)$.
		The functional equation is
		\[
		\widehat Z_{i1}(s,\nu,\ell)
		=|\nu|^{2-2s}\overline{\alpha(\nu)}^{\,\ell}
		\sum_{i'} A_{i'i}(2^{-s},\ell)\widehat Z_{i'1}(2-s,\nu,-\ell).
		\]
		Here $A_{i'i}$ denotes the corresponding rational matrix entry in
		Eq.~(4.28) of that paper.
		Proposition~4.3 shows that at angular weight $1$, $\widehat Z_{i1}$
		is entire and bounded in vertical strips.
		
		Move the Mellin contour from $\Re z=9/8$ to $\Re z=-1/8$.
		The growth is polynomial on both boundaries, by absolute convergence
		on the right and by the functional equation and Stirling's formula on
		the left.
		With the archimedean factor removed, the growth inside the strip is
		at most exponential, so the Phragm\'en--Lindel\"of principle gives
		the polynomial bound needed for the contour shift.
		After reflection, with $t=1-z$, the conjugate of the angular factor is
		$\alpha(\nu)=\alpha(r)^3$, and
		\[
		X^z |\nu|^{1-2z}
		=X |\nu|^{-1}\left(\frac{|\nu|^2}{X}\right)^t,
		\qquad |\nu|^2=\frac{|r|^6}{4}.
		\]
		
		To specify the terms in \eqref{eq:reflection}, use
		\cite[Remark~4.2]{ddhl} to write the series at each reflected cusp as
		\[
		|\lambda|^8\,\psi^{(4)}_{i'1}(s,r^3/2,-1)
		=\sum_{j:\,i'(j)=i'} \overline{\eta_j(r)}\,2^{-m_j s}
		\Psi_{\beta_j}(s,\kappa_j r^3,-1).
		\]
		Here $i'(j)$ is the cusp index for the $j$th term,
		$m_j$ is its dyadic exponent, and $\overline{\eta_j(r)}$ is the
		corresponding finite coefficient in Eq.~(4.15) or Eq.~(4.16).
		With $\nu=r^3/2$, the additive character in Eq.~(4.15)
		depends only on $r\bmod\lambda^4$.
		By Eqs.~(4.20)--(4.21) and $v_\lambda(\nu)=-2$, it suffices to
		consider the $(a,m_j)$ terms in Eq.~(4.16) with $0\le a\le3$ and
		$2\le m_j\le3$.
		The additive factor in Eq.~(4.17) that varies with $r$ reduces to
		\[
		e_0\left(\frac{h r^3}{i^a\lambda^{m_j+4}}\right),
		\]
		Here $h$ ranges over the residue classes specified in that equation.
		This factor depends only on $r\bmod\lambda^{m_j+4}$.
		Thus $\eta_j(r)$ is bounded on odd $r$ and depends only on
		$r\bmod\lambda^7$, and it can be set to zero on even $r$.
		Put
		\[
		B_j(t)=-2J(\beta)\,4^{-t}2^{-m_j(t+1/2)}
		\overline{A_{i'(j)i_\beta}(2^{\bar t-3/2},1)}
		\frac{G_\infty(t+1/2,-1)}{G_\infty(3/2-t,1)} .
		\]
		The factors $2$ and $4^{-t}$ come from $|\nu|^2=|r|^6/4$.
		For positive real $v,y$, define the kernel by
		\begin{equation}\label{eq:kernel-definition}
			\mathcal W_{j,b}(v,y)=\frac{1}{2\pi i}\int_{(a)}
			\widehat w_{v,b}(1-t)B_j(t)y^{-t}\,\mathrm dt,\qquad a>1.
		\end{equation}
		Expanding the reflected factor $L_\lambda(4t-1,\alpha^{-4})$ as a sum over
		$e$ gives \eqref{eq:reflection}, including the weight
		$|e|^2\,\overline{\alpha(e)}^{\,4}$.
		
		Before reflection, the possible poles of the rational matrix
		coefficients lie on the lines with real parts $1,5/4$.
		Thus $B_j(t)$ is holomorphic in $\Re t>1$ and has polynomial growth
		on fixed vertical lines.
		Integrating by parts in $x$ in the formula for $\widehat w_{v,b}(1-t)$
		gives decay of any polynomial order on these lines while retaining the
		factor $(1+v/A_b)^{-L}$.
		Euler derivatives in $v$ and $y$ preserve this property.
		The line of integration in \eqref{eq:kernel-definition} can therefore be moved
		within $\Re t>1$.
		For $0<y\le1$, take $a=9/8$, and for $y\ge1$, take $a$ to be a
		power of $2$ greater than $\max(1,B)$, giving \eqref{eq:kernel}.
	\end{proof}
	
	\subsection{Second Poisson summation and the upper bound}
	
	\begin{lemma}[Nonzero angular weight and the zero-frequency term]\label{lem:angular}
		Let $R\ge1$ be real and let $h$ be a smooth function supported
		on $[R/2,2R]$.
		Let $B_h$ be a common upper bound for $|(v\partial_v)^a h(v)|$
		for $a=0,\ldots,8$.
		Let $P$ be a function on the Gaussian integers, periodic modulo
		$b\lambda^7$, with $|P|\le1$.
		For any $\kappa_j$ in \eqref{eq:reflection} and any odd primary $n$,
		we have
		\begin{equation}\label{eq:angular}
			\left|\sum_{r\in\mathbb Z[i]\setminus\{0\}}P(r)\alpha(r)^3
			\frac{\overline{g_4(\kappa_jr^3,n)}}{|n|}h(|r|^2)\right|
			\ll_\varepsilon |b|^2\,|n|^{1+2\varepsilon}B_h.
		\end{equation}
	\end{lemma}
	
	\begin{proof}
		Expand the Gauss sum using its multiplication formula at each
		$\tau^k\Vert n$, and write $q=|\tau|^2$.
		Since $2\kappa_j$ is coprime to $\tau$, the local factor has the same
		form as the first Gauss sum in Lemma~\ref{lem:prime-power}.
		The normalized expression has at most $2^{\omega(n)}$ components.
		For each component, $a_0$ is the primary Gaussian integer determined
		by the local divisibility conditions.
		On lattice points of the form $r=a_0t$, where $t$ is a Gaussian
		integer, the component is a primitive character $\psi(t)$ modulo
		a squarefree Gaussian integer $f$, multiplied by a constant $C$.
		The character $\psi$ is zero on nonunits modulo $f$ and identically
		$1$ when $f=1$.
		A nontrivial local character has a constant factor of modulus
		$q^{(k-1)/2}$ and contributes a factor of $\tau$ to $f$.
		A constant component has modulus at most $q^{k/2}$ and contributes
		no factor to $f$.
		Therefore $|C|\,|f|\le|n|$.
		
		Apply Poisson summation in $t$ to this component.
		For the finite part $P(a_0t)\psi(t)$, we can use
		$q_0=\operatorname{lcm}(b\lambda^7,f)$ as the period modulus, taking
		any generator of this least common multiple ideal.
		Separate the part of $f$ coprime to $b\lambda^7$ using the Chinese
		remainder theorem, and apply the primitive Gauss sum formula to bound
		the finite Fourier sum by $|b\lambda^7|^2|f|$.
		The continuous part is a radial weight on the norm scale $R/|a_0|^2$
		multiplied by $\alpha(t)^3$ and by the constant
		$\alpha(a_0)^3$.
		The integral at zero frequency vanishes because
		$\int_0^{2\pi}e^{3i\vartheta}\,\mathrm d\vartheta=0$.
		For the remaining frequencies, integration by parts shows that the
		sum of the absolute values of the Fourier weights, divided by $|q_0|^2$,
		is $O(B_h)$.
		This uses
		\[
		\frac{R}{|a_0|^2\,|q_0|^2}\sum_{m\in\mathbb Z[i]\setminus\{0\}}
		\left(1+|m|\sqrt{\frac{R}{|a_0|^2\,|q_0|^2}}\right)^{-4}\ll 1.
		\]
		Summing over dyadic ranges of $|m|$ gives this bound, and the
		constant is independent of the ratio of scales.
		Finally, multiply by $|C|$, sum over the components, and use
		$2^{\omega(n)}\ll_\varepsilon|n|^{2\varepsilon}$ to obtain
		\eqref{eq:angular}.
	\end{proof}
	
	\begin{proof}[Proof of Proposition~\ref{prop:dual}]
		In \eqref{eq:poisson}, take a smooth dyadic partition of unity
		$\sum_{R=2^j,\ j\ge0}\vartheta_R(v)=1$ ($v\ge1$), where $R$
		runs over the nonnegative integer powers of $2$.
		Each $\vartheta_R$ is supported on $[R/2,2R]$, and its Euler
		derivatives of every order are uniformly bounded.
		Substitution of \eqref{eq:reflection} gives the factor
		$\eta_j(r)\widehat J(r)g_-(r,b)/|b|^2$, which is bounded and periodic
		modulo $b\lambda^7$ and vanishes on even $r$.
		Absorb $(|r|^2/R)^{-3/2}$ into the smooth weight and apply
		Lemma~\ref{lem:angular}.
		In \eqref{eq:kernel}, take $B=L=8$, and in what follows it suffices
		to consider $0<\varepsilon\le1/8$.
		For $0<y\le1$, set $\Phi(y)=y^{-9/8}$, and for $y>1$, set
		$\Phi(y)=y^{-8}$.
		By \eqref{eq:kernel},
		\begin{equation}\label{eq:dyadic-r}
			\begin{aligned}
				&\left|\sum_{\lambda\nmid r}\vartheta_R(|r|^2)\widehat J(r)
				\frac{g_-(r,b)}{|b|^2}\mathcal S_\beta(r;b)\right|\\
				&\qquad\ll_\varepsilon |b|^2\,X R^{-3/2}(1+R/A_b)^{-8}
				\sum_{n,e} |e|^2\,|n|^{1+2\varepsilon}
				\Phi\left(\frac{|n|^2\,|e|^8X}{R^3}\right),
			\end{aligned}
		\end{equation}
		The radial derivatives in $r$ are controlled by \eqref{eq:kernel},
		since the second argument of the kernel is a constant times $|r|^{-6}$.
		
		Let $D_0$ run over the nonnegative integer powers of $2$.
		On the range $|n|^2\,|e|^8\sim D_0$, the total weight in
		\eqref{eq:dyadic-r} with the kernel omitted is
		\[
		\ll_\varepsilon
		\sum_e |e|^2\left(\frac{D_0}{|e|^8}\right)^{3/2+\varepsilon}
		\ll_\varepsilon D_0^{3/2+\varepsilon}.
		\]
		This uses the fact that the number of primary Gaussian integers of
		norm at most $Z$ is $O(Z)$.
		Since $9/8<3/2+\varepsilon<8$, summing over dyadic $D_0$ gives
		\[
		O_\varepsilon\bigl((R^3/X)^{3/2+\varepsilon}\bigr).
		\]
		When $R^3<X$, the large-argument bound for $\Phi$ gives the same
		estimate directly.
		Since $3+3\varepsilon<8$, summing over $R$ shows that the inner sum
		in \eqref{eq:poisson} for each $b$ is
		\[
		O_\varepsilon\left((XHY)^\varepsilon |b|^2\,X^{-1/2}A_b^3\right).
		\]
		Here the power losses are absorbed into $(XHY)^\varepsilon$.
		Since $|b|^2\asymp Y$ and $A_b=X|b|^2/H$, we have
		\[
		\sum_b \bigl|W_1(|b|^2/Y)\bigr|\,|b|^2\,A_b^3\ll (X/H)^3Y^5 .
		\]
		Multiplication by the prefactor $H/(|\lambda|^8\sqrt X)$ gives
		\eqref{eq:dual-bound}.
	\end{proof}
	
	\section{Proof of the estimate for the nonprincipal contribution}\label{sec:nonprincipal}
	
	\subsection{Double Mellin inversion}
	
	For primary $u$ and complex numbers $z,w$ satisfying $\Re z,\Re w>1$,
	define an Euler product $P_u(z,w)$.
	Its factor at an odd Gaussian prime $\tau$ is
	\[
	\frac{1-\overline{\xi_u(\tau)}q^{-w}}{1-\xi_u(\tau)q^{-z}}
	\sum_{k,\ell\ge0}F_{\tau^k,\tau^\ell}(u)q^{-kz-\ell w},
	\qquad q=|\tau|^2.
	\]
	
	\begin{lemma}\label{lem:euler}
		In $\Re z,\Re w>1$,
		\begin{equation}\label{eq:double}
			\sum_{c,b}\frac{F_{c,b}(u)}{|c|^{2z}|b|^{2w}}
			=\frac{L_\lambda(w,\overline{\xi_u})}{L_\lambda(z,\xi_u)}P_u(z,w).
		\end{equation}
		The function $P_u$ is holomorphic in the region defined by
		$\Re z,\Re w>0$ and $\Re(z+w)>1$.
		If $\sigma=\Re z\ge1/2$, $\theta=\Re w\ge1/2$, and $\sigma+\theta\ge5/4$,
		then
		\begin{equation}\label{eq:correction}
			\begin{gathered}
				|P_u(z,w)|\ll_\varepsilon |u|^{2(\varepsilon+\Lambda(\sigma,\theta))},\\
				\Lambda(\sigma,\theta)=\max(0,1-\tfrac32\sigma,1-\sigma-\tfrac13\theta).
			\end{gathered}
		\end{equation}
	\end{lemma}
	
	\begin{proof}
		By \eqref{eq:crt}, the left-hand side factors as an Euler product.
		If $\tau\nmid u$, Lemma~\ref{lem:prime-power} gives the local factor
		\[
		(1-\overline{\xi_u(\tau)}q^{-w})^{-1}-\xi_u(\tau)q^{-z},
		\qquad q=|\tau|^2.
		\]
		Factoring out the two local $L_\lambda$ factors leaves
		$1+q^{-z-w}/(1-\xi_u(\tau)q^{-z})$.
		The product of these factors converges absolutely and locally
		uniformly in the stated region.
		If $\tau\mid u$, the $q^v$ bound in \eqref{eq:local} gives absolute
		convergence of the local series in $\Re z,\Re w>0$.
		The finite correction from factoring out the $L_\lambda$ factors
		is also holomorphic in this region.
		
		To prove \eqref{eq:correction}, write $v=v_\tau(u)$ and put
		$a=\min(v,k,\lfloor2k/3\rfloor+\ell)$.
		From $k\ge a$ and $\ell\ge a-2k/3$, we obtain
		\[
		k\sigma+\ell\theta
		\ge a\min(3\sigma/2,\sigma+\theta/3).
		\]
		When $k\ge3a/2$, it suffices to use $k\sigma$.
		When $a\le k\le3a/2$, evaluate the linear lower bound in $k$ at the
		two endpoints.
		Thus the terms with $k+\ell\le4(v+1)$ are all
		$O(q^{v\Lambda(\sigma,\theta)})$, and there are $O((v+1)^2)$ such terms.
		The other terms are controlled by the geometric decay of
		$q^{v-(k+\ell)/2}$.
		The finite Euler corrections are uniformly bounded in the region above.
		A divisor bound absorbs $(v+1)^2$ and the constant at each bad
		prime divisor, giving \eqref{eq:correction}.
	\end{proof}
	
	Applying double Mellin inversion to \eqref{eq:Bd} and substituting
	\eqref{eq:double} gives
	\begin{equation}\label{eq:row-integral}
		\begin{aligned}
			\mathcal B_d(u)
			&=\frac{1}{(2\pi i)^2}
			\int_{(1+\varepsilon)}\int_{(1+\varepsilon)}
			\widehat W(z)\widehat W_1(w) \frac{X^zY^w}{|d|^{8z}}\\
			&\hspace{2.5em}\times
			\frac{L_\lambda(w,\overline{\xi_u})}{L_\lambda(z,\xi_u)}
			P_u(z,w)\,\mathrm dw\,\mathrm dz .
		\end{aligned}
	\end{equation}
	We also use this formula in Section~\ref{sec:principal} to compute the
	principal character contribution.
	
	\subsection{Reciprocal bounds and zero density}
	
	\begin{lemma}[Reciprocal bound in a zero-free rectangle]\label{lem:reciprocal}
		Fix a real parameter $0<t<1$, put $T=X^t$, and let $\xi$ be a nonprincipal
		primitive Hecke character of finite order with conductor norm $O(X)$.
		If $L(s,\xi)$ has no zeros in the region $\Re s>a,\ |\Im s|\le2T+4$, where
		$5/8\le a\le1$, then for any fixed $\Delta>0,\varepsilon>0$,
		\begin{equation}\label{eq:reciprocal}
			L_\lambda(s,\xi)^{-1}\ll_{t,\Delta,\varepsilon} X^\varepsilon
			\qquad(\Re s\ge a+\Delta,\ |\Im s|\le T).
		\end{equation}
	\end{lemma}
	
	\begin{proof}
		The Euler series in $\Re s>1$ fixes a branch of $\log L$, which we
		extend to the given zero-free rectangle.
		The standard Hecke functional equation and absolute convergence give
		polynomial bounds in fixed vertical strips, so $\log|L|\ll\log(2X)$
		throughout the rectangle.
		For real $v$ with $|v|\le3T/2$, apply the Borel--Carath\'eodory
		inequality on disks centered at $2+iv$.
		Taking radii $2-a-\Delta/4$ and $2-a-\Delta/2$ gives
		$\log L(s)\ll_\Delta\log(2X)$ for
		$\Re s\ge a+\Delta/2,\ |\Im s|\le3T/2$.
		Farther to the right, absolute convergence applies directly.
		Here it suffices to consider $\Delta<1/4$.
		
		For $\Re s\ge1+\Delta/4$, the Euler series gives $\log L=O_\Delta(1)$.
		If $a+\Delta<1+\Delta/4$, then for each $|v|\le T$, apply the maximum
		modulus principle to the following expression in the rectangle with
		real part between $a+\Delta/2$ and $1+\Delta/4$ and imaginary part
		between $-3T/2$ and $3T/2$.
		\[
		\log L(s)\exp((s-iv)^2)
		\bigl(\log(2X)\bigr)^{-(1+\Delta/4-s)/(1-a-\Delta/4)}.
		\]
		This expression is $O_\Delta(1)$ on both vertical sides.
		On the horizontal sides, $|\Im s-v|\ge T/2$, so it is
		$O_\Delta(\log(2X)\exp(-T^2/4))=O_{t,\Delta}(1)$.
		Therefore, for $\Re s\ge a+\Delta$ and $\Im s=v$, we have
		$|\log L(s)|\ll_{t,\Delta}(\log(2X))^{1-\Delta/2}$.
		Exponentiation gives $L^{-1}\ll_{t,\Delta,\varepsilon}X^\varepsilon$.
		The omitted local factor and its reciprocal are uniformly bounded
		in this region.
	\end{proof}
	
	For $1/2\le\sigma\le1$, write $q_0=2\sigma-1$ and let $E(\sigma)$
	be the norm exponent in \cite[Corollary~1.5]{defaveri} for $n=4$.
	The explicit formula is
	\[
	E(\sigma)=
	\begin{cases}
		1-\dfrac{q_0}{4(1-q_0)},& 1/2\le\sigma\le3/4,\\[4pt]
		1-\dfrac{q_0(4q_0-1)}{24q_0-16q_0^2-6},
		&3/4\le\sigma\le7/8,\\[4pt]
		\dfrac{6(1-\sigma)}{1+4(1-\sigma)},&7/8\le\sigma\le1 .
	\end{cases}
	\]
	Given $Q>1$, $T\ge1$, and $5/8\le\sigma<1$, consider fourth-power-free
	primary Gaussian integers $u_0$ with $|u_0|^2\sim Q$.
	The number of such $u_0$ for which $L(s,\xi_{u_0})$ has at least one
	zero in the rectangle $\sigma\le\Re s\le1$, $|\Im s|\le 2T+4$
	is at most
	\begin{equation}\label{eq:zero-density}
		O_\varepsilon\bigl(Q^{E(\sigma)+\varepsilon}T^{2+\varepsilon}\bigr).
	\end{equation}
	Example~2.3 of that paper takes $S=\{(\lambda)\}$ and gives the
	character indexed by $(u_0)$ as $\left(\frac{u_0}{a}\right)_4$, with associated
	primitive character $\xi_{u_0}$.
	The sum in Corollary~1.5 ranges over fourth-power-free ideals coprime
	to $S$, so it includes the ideals $(u_0)$ here.
	Substituting $n=4$ and field degree $2$ gives the norm exponent
	$E(\sigma)$ above and the height exponent
	$1+2(1-\sigma)/(3-2\sigma)+\varepsilon\le2+\varepsilon$.
	The total number of zeros up to height $2T+4$ bounds the number
	of such primary generators.
	The omitted finite Euler factors are nonzero in $\Re s>0$, so the
	same estimate applies to the primitive $L$-functions.
	
	\begin{lemma}[Comparison of zero-density exponents]\label{lem:density}
		Let $5/8\le\sigma\le599/600$.
		Take $E=E(\sigma)$, with the choice $E=1$ also allowed at $\sigma=5/8$.
		Then
		\begin{equation}\label{eq:density-margin}
			\sigma+\frac{65}{144}\bigl(\max(E-\tfrac14,0)+\Lambda(\sigma,1)\bigr)
			\le \frac{599}{600}.
		\end{equation}
	\end{lemma}
	
	\begin{proof}
		By \eqref{eq:correction}, $\Lambda(\sigma,1)=\max(0,1-\tfrac32\sigma)$.
		When $\sigma=5/8$ and $E=1$, the left-hand side is $2285/2304<599/600$.
		If $E=E(\sigma)<1/4$, then $\sigma>2/3$, in which case the left-hand
		side equals $\sigma$.
		In the remaining cases it suffices to bound
		\[
		B(\sigma)=\sigma+\frac{65}{144}\bigl(E(\sigma)+\Lambda(\sigma,1)-\tfrac14\bigr).
		\]
		On $[5/8,2/3]$, we have $E(\sigma)\le11/12$ and
		$\sigma+(65/144)\Lambda(\sigma,1)\le2/3$, so $B(\sigma)\le209/216$.
		On $[2/3,3/4]$, we have $B'(\sigma)=1-(65/144)/(8(1-\sigma)^2)>0$,
		so $B(\sigma)\le B(3/4)=281/288$.
		
		On $[3/4,7/8]$, put $x=q_0-33/50$, so that $x\ge-4/25$.
		Substituting the expression for $E(\sigma)$ gives
		\begin{align*}
			&(24q_0-16q_0^2-6)\left(\frac{599}{600}-B(\sigma)\right)\\
			&\qquad=8x^3+\frac{139}{45}x^2+\frac{589}{90000}x+\frac{13}{750000}\\
			&\qquad\ge\frac95x^2-\frac{|x|}{150}+\frac1{60000}>0.
		\end{align*}
		Here $24q_0-16q_0^2-6\ge2$, and the last quadratic in $|x|$ has
		negative discriminant.
		On $[7/8,1]$, we have $B'(\sigma)=1-6(65/144)/(1+4(1-\sigma))^2<0$,
		so $B(\sigma)\le B(7/8)=569/576$.
	\end{proof}
	
	\begin{proof}[Proof of Proposition~\ref{prop:nonprincipal}]
		Let $Q$ be a positive real number satisfying $|u_0|^2\sim Q$.
		Take sufficiently small positive real numbers $t,\Delta$ and put
		$T=X^t$.
		Consider only zeros of $L(s,\xi_{u_0})$ with $|\Im s|\le2T+4$.
		Put the characters for which the largest real part of a zero is at
		least $599/600$ in a separate group and take $\sigma=599/600$.
		Assign the remaining characters for which this maximum is at least
		$5/8$ to intervals of width $\Delta$ according to that maximum,
		and take $\sigma$ to be the left endpoint of each interval.
		For all the groups above, take $E=E(\sigma)$.
		For the group with no zero of real part at least $5/8$, take
		$\sigma=5/8$ and $E=1$.
		By \eqref{eq:zero-density} or by counting lattice points, each group
		contains at most $O_\varepsilon(X^\varepsilon Q^E)$ choices of $u_0$.
		Each $u_0$ corresponds to $O(H^{1/4}|d|^{-2}Q^{-1/4})$ choices of the
		primary $f$.
		Therefore, the number of primary Gaussian integers $u$ in each
		group satisfying $|u_0|^2\sim Q$ and $|u|^2\,|d|^8\le2H$ is
		\begin{equation}\label{eq:density-count}
			O_\varepsilon\left(\frac{X^\varepsilon H^{1/4}Q^{E-1/4}}{|d|^2}\right).
		\end{equation}
		
		For the group with $\sigma=599/600$, keep the $z$ contour in
		\eqref{eq:row-integral} fixed and move the $w$ contour to $\Re w=13/25$.
		For a nonprincipal character, the numerator is holomorphic, and
		$\Lambda=0$ in \eqref{eq:correction}.
		On $\Re w=-\varepsilon$, the Hecke functional equation gives a bound
		of $Q^{1/2+\varepsilon}$ times a polynomial in the height.
		Interpolating with the absolute convergence bound on
		$\Re w=1+\varepsilon$ gives
		\[
		L_\lambda(w,\overline{\xi_u})
		\ll_\varepsilon Q^{6/25+\varepsilon}(1+|\Im w|)^B
		\qquad(\Re w=13/25),
		\]
		Here $B$ is a constant.
		In this case $E(599/600)=3/302$, so $E+6/25-1/4=-1/15100$.
		Summing over $Q$ and $d$ using \eqref{eq:density-count}, we find that
		the contribution of this group is
		\[
		\ll_\varepsilon X^{1+\varepsilon}Y^{13/25}H^{1/4}\sum_d|d|^{-4}
		\ll_\varepsilon X^{599/600+\varepsilon}YH^{1/4}.
		\]
		
		For the remaining groups, keep $\Re w=1+\varepsilon$, truncate the
		$z$ integral at height $T$, and shift the contour to $\Re z=\sigma+2\Delta$.
		Taking $\Delta<1/1200$ puts the new line to the left of $1$.
		In Lemma~\ref{lem:reciprocal}, take $a=\sigma+\Delta$.
		The resulting reciprocal bound applies on the new line and on the
		horizontal segments with $|\Im z|=T$ and $\Re z\ge\sigma+2\Delta$.
		The tails removed by truncation lie on the original line $\Re z=1+\varepsilon$,
		where $L_\lambda(z,\xi_u)^{-1}$ is bounded uniformly by the absolutely
		convergent Euler product.
		The other factors grow at most polynomially.
		Because $\widehat W$ decays rapidly, the horizontal segments and the
		tails each contribute $O(X^{-C})$ for arbitrary $C$.
		
		Using \eqref{eq:correction}, we absorb $(X/|d|^8)^{2\Delta}$ into
		$X^\varepsilon$ and obtain the following bound for each integral indexed by $(d,u)$.
		\[
		O_\varepsilon\left(\frac{X^{\sigma+\varepsilon}Y H^{\Lambda(\sigma,1)}}{|d|^{8(\sigma+\Lambda(\sigma,1))}}\right).
		\]
		Write $A=\max(E-\tfrac14,0)+\Lambda(\sigma,1)$.
		Using $Q\ll H/|d|^8$ in \eqref{eq:density-count} and then
		summing over dyadic $Q$, we obtain the following bound for the
		contribution of this group.
		\[
		\ll_\varepsilon X^\varepsilon Y X^\sigma H^{1/4+A}
		\sum_d |d|^{4-8(\sigma+A)}.
		\]
		When $\sigma\le3/4$, we have $E\ge3/4$, so $\sigma+A>3/4$.
		When $\sigma>3/4$, the bound $A\ge0$ gives the same inequality.
		Thus the sum over $d$ converges.
		By Lemma~\ref{lem:density}, $\sigma+(65/144)A\le599/600$, so
		the expression above is $O_\varepsilon(X^{599/600+\varepsilon}YH^{1/4})$.
		Choose $t,\Delta$ and all exponent losses sufficiently small, then
		sum over the groups to obtain the desired estimate.
	\end{proof}
	
	\section{Derivation of the main term}\label{sec:principal}
	
	\subsection{The averaged Euler product}
	
	When $u=f^4$ and $v_\tau(f)=j$, the quotient $f/\tau^j$ is coprime to
	$\tau$.
	By the scaling formula in Lemma~\ref{lem:crt},
	$F_{\tau^k,\tau^\ell}(f^4)=F_{\tau^k,\tau^\ell}(\tau^{4j})$, since
	$(f/\tau^j)^4$ is a fourth power of a unit modulo $\tau^{k+\ell}$.
	Thus $D_{\tau,j}(z)$ in \eqref{eq:local-D} is the local factor on the
	left-hand side of \eqref{eq:double} when $u=f^4$ and $w=1$.
	Factoring out the $L_\lambda$ factors for the principal character gives
	\begin{equation}\label{eq:principal-product}
		P_{f^4}(z,1)=\prod_{\tau\nmid\lambda}P_{\tau,v_\tau(f)}(z)
		\qquad(\Re z>1).
	\end{equation}
	
	\begin{lemma}[Nonvanishing of the averaged Euler factors]\label{lem:nonzero}
		The Euler product $\mathcal Q(z)$ is holomorphic and nonzero in
		$\Re z>1-2^{-6}$.
	\end{lemma}
	
	\begin{proof}
		Fix an odd Gaussian prime $\tau$ and write $q=|\tau|^2$.
		If $j=0$, Lemma~\ref{lem:prime-power} gives
		$D_{\tau,0}(z)=(1-q^{-1})^{-1}-q^{-z}$, so
		\[
		P_{\tau,0}(z)=1+\frac{q^{-z-1}}{1-q^{-z}}.
		\]
		If $j\ge1$, the local factor consists of the constant term, the term
		with $(k,\ell)=(1,1)$, and the terms with $k\ge4$.
		Write $\sigma=\Re z\ge1-2^{-6}$.
		The absolute value of the $(1,1)$ term is at most
		$(q-1)q^{-\sigma-1}$.
		For $k\ge4$, put $m=\lceil k/3\rceil$.
		By \eqref{eq:local},
		\begin{align*}
			\sum_{\ell\ge0}|F_{\tau^k,\tau^\ell}(\tau^{4j})|q^{-k\sigma-\ell}
			&\le 4q^{-k\sigma}
			\left(\sum_{0\le\ell<m}q^{k-m}+\sum_{\ell\ge m}q^{k-\ell}\right)\\
			&=4q^{-m+(1-\sigma)k}\left(m+1+\frac1{q-1}\right).
		\end{align*}
		Since $k\le3m$ and each $m$ corresponds to at most three values of
		$k$, the total contribution from $k\ge4$ is at most
		\begin{equation}\label{eq:euler-tail}
			12\sum_{m\ge2}\left(m+1+\frac1{q-1}\right)q^{-61m/64}
			\ll q^{-61/32}.
		\end{equation}
		It follows that $P_{\tau,j}(z)=1+O(q^{-\min(\sigma,1)})$ uniformly
		for $j\ge1$.
		The local average in \eqref{eq:mean-euler} is therefore
		$1+O(q^{-1-\min(\sigma,1)})$, and the product converges absolutely
		and locally uniformly in $\sigma>1-2^{-6}$.
		
		To prove that each local average is nonzero, it suffices to show
		\begin{equation}\label{eq:euler-close}
			\left|(1-q^{-1})\sum_{j\ge0}q^{-j}D_{\tau,j}(z)-1\right|<1.
		\end{equation}
		The term with $j=0$, together with the constant terms for $j\ge1$,
		differs from $1$ by at most $q^{-1}+q^{-\sigma}$.
		The contribution of the $(1,1)$ term is at most $q^{-\sigma-1}$,
		and the terms covered by \eqref{eq:euler-tail} contribute at most
		\[
		12q^{-1}\sum_{m\ge2}\left(m+1+\frac1{q-1}\right)q^{-61m/64}.
		\]
		The norm of an odd Gaussian prime satisfies $q\ge5$.
		Putting $x=q^{-61/64}$ gives $q^{-\sigma}\le x<7/32$, since
		$(9/8)^{16}>5$ and $(9/8)^{3/4}\le35/32$.
		The sum of the first two bounds is at most $1/5+(6/5)(7/32)<15/32$.
		Using the bound $x\le7/32$ in
		\[
		\sum_{m\ge2}\left(m+1+\frac1{q-1}\right)x^m
		=x^2\left(\frac{3+1/(q-1)}{1-x}+\frac{x}{(1-x)^2}\right)
		\]
		shows that the last contribution is less than $17/32$.
		Thus \eqref{eq:euler-close} holds.
		Multiplying back by $(1-q^{-1})/(1-q^{-z})$ preserves nonvanishing.
		By absolute convergence, $\mathcal Q$ is holomorphic and nonzero.
	\end{proof}
	
	\subsection{Lattice averaging and the residue}
	
	\begin{lemma}[Lattice average over fourth powers]\label{lem:lattice}
		We have $\operatorname*{Res}_{w=1}\zeta_\lambda(w)=\pi/8$.
		For a positive real $R$ and a complex number $z$ satisfying
		$\Re z=1+\varepsilon$, we have
		\begin{equation}\label{eq:lattice}
			\sum_f V\bigl((|f|^2/R)^4\bigr)P_{f^4}(z,1)
			=\frac{\pi R}{8}\mathcal Q(z)\int_0^\infty V(x^4)\,\mathrm dx
			+O_\varepsilon(R^{1/2}),
		\end{equation}
		Here the error is uniform in $\Im z$.
	\end{lemma}
	
	\begin{proof}
		The primary Gaussian integers lie in one lattice coset modulo
		$\lambda^3$, with fundamental area $8$.
		Estimating the lattice cells that meet the boundary shows that
		the number of primary Gaussian integers of norm at most $R$ is
		$\pi R/8+O(R^{1/2})$.
		The same coefficient is also the residue of $\zeta_\lambda$ at $1$,
		as follows by partial summation.
		Thus for any $R>0$, summing over primary $g$ gives
		\begin{equation}\label{eq:lattice-basic}
			\sum_g V\bigl((|g|^2/R)^4\bigr)
			=\frac{\pi R}{8}\int_0^\infty V(x^4)\,\mathrm dx+O(R^{1/2}).
		\end{equation}
		When $R<1$, the support of $V$ also gives the same error bound directly.
		
		The proof of Lemma~\ref{lem:nonzero} shows that on $\Re z=1+\varepsilon$,
		the product $\prod_{\tau\nmid\lambda}P_{\tau,0}(z)$ is bounded and
		$P_{\tau,0}(z)$ is uniformly bounded away from zero.
		Define a multiplicative function $a_z$ on the primary Gaussian
		integers by $a_z(1)=1$ and
		\[
		a_z(\tau^j)=\frac{P_{\tau,j}(z)-P_{\tau,j-1}(z)}{P_{\tau,0}(z)}
		\quad(j\ge1).
		\]
		By telescoping the local factors and using $a_z(\tau^j)\ll|\tau|^{-2}$,
		we have
		\[
		P_{f^4}(z,1)=\left(\prod_{\tau\nmid\lambda}P_{\tau,0}(z)\right)\sum_{a\mid f}a_z(a),
		\qquad \sum_a |a_z(a)|\,|a|^{-1}\ll1 .
		\]
		For each $a$, write $f=ag$ and apply \eqref{eq:lattice-basic} at scale
		$R/|a|^2$.
		The errors can be summed absolutely, giving $O(R^{1/2})$.
		The coefficient of the main term is
		$\left(\prod_{\tau\nmid\lambda}P_{\tau,0}(z)\right)\sum_a a_z(a)/|a|^2$.
		For each odd Gaussian prime $\tau$ of norm $q=|\tau|^2$, we have
		\[
		P_{\tau,0}(z)\left(1+\sum_{j\ge1}q^{-j}a_z(\tau^j)\right)
		=(1-q^{-1})\sum_{j\ge0}q^{-j}P_{\tau,j}(z).
		\]
		Thus the coefficient of the main term is $\mathcal Q(z)$, proving
		\eqref{eq:lattice}.
	\end{proof}
	
	\begin{proof}[Proof of Proposition~\ref{prop:principal}]
		In \eqref{eq:row-integral} for the principal character, keep
		$\Re z=1+\varepsilon$ and move the $w$ contour to $\Re w=13/25$,
		crossing the simple pole of $\zeta_\lambda(w)$ at $w=1$.
		Estimate the integral on the new line using the polynomial growth of
		$\zeta_\lambda(w)$ and \eqref{eq:correction} with $\Lambda=0$.
		There are $O(H^{1/4}/|d|^2)$ choices of $f$, and summing over $d$ gives
		$O_\varepsilon(X^{1+\varepsilon}Y^{13/25}H^{1/4})$.
		
		The residue term is
		\[
		\frac{\pi}{8}Y\widehat W_1(1)\frac1{2\pi i}
		\int_{(1+\varepsilon)}
		\frac{\widehat W(z)X^z}{\zeta_\lambda(z)}
		\sum_d\frac{|d|^6\alpha(d)^4}{|d|^{8z}}
		\sum_f V\left(\frac{|f|^8\,|d|^8}{H}\right)
		P_{f^4}(z,1)\,\mathrm dz .
		\]
		Apply \eqref{eq:lattice} to the sum over $f$ at scale $H^{1/4}/|d|^2$.
		The error is controlled by the convergence of
		$\sum_d|d|^{-3-8\varepsilon}$, so it is
		$O_\varepsilon(X^{1+\varepsilon}YH^{1/8})$.
		In the main term, the weight of $d$ becomes
		$\alpha(d)^4/|d|^{8z-4}$, and summing these weights gives
		$L_\lambda(4z-2,\alpha^4)$.
		The remaining factors combine to give $K_0YH^{1/4}$, yielding
		\eqref{eq:principal}.
	\end{proof}
	
	\section{Proof of the main theorem}\label{sec:conclusion}
	
	\begin{proof}[Proof of Theorem~\ref{thm:main}]
		Use the decomposition \eqref{eq:I-parts} to combine
		Propositions~\ref{prop:dual}, \ref{prop:nonprincipal}, and
		\ref{prop:principal}.
		By \eqref{eq:parameters},
		\[
		\frac{X^2Y^5H^{-2}}{XYH^{1/4}}=X^{-1/576},\qquad
		Y^{-12/25}=X^{-1/600},\qquad H^{-1/8}=X^{-65/1152}.
		\]
		In Proposition~\ref{prop:dual}, take the loss parameter to be
		$\varepsilon/4$, and use $H,Y\le X$ to obtain
		$(XHY)^{\varepsilon/4}\le X^\varepsilon$.
		Combining \eqref{eq:nonprincipal} and \eqref{eq:principal}, we find
		that after division by $XYH^{1/4}$ all the errors are
		$O_\varepsilon(X^{-1/600+\varepsilon})$.
		Since $K_0>0$, we obtain
		\begin{equation}\label{eq:s-saving}
			S_W(X)\ll_{\varepsilon,W} X^{599/600+\varepsilon}\qquad(X\ge2).
		\end{equation}
		
		Since \eqref{eq:s-saving} holds for every $\varepsilon>0$ and
		$S_W(X)=0$ ($0<X<1/2$), the right-hand side of
		\eqref{eq:mellin-identity} converges absolutely and locally uniformly
		in $\Re z>599/600$, giving a holomorphic continuation of
		$\widehat W(z)D(z)$.
		
		By Lemma~\ref{lem:nonzero}, $\mathcal Q(z)$ is holomorphic and
		nonzero in $\Re z>1-2^{-6}$.
		In the same region, $4\Re z-2>1$, so $L_\lambda(4z-2,\alpha^4)$
		is nonzero by its absolutely convergent Euler product.
		For any $z_0$ with $\Re z_0>599/600$, take a nonnegative and nonzero
		$\varphi\in C_c^\infty((1,2))$ and put $W(x)=x^{-z_0}\varphi(x)$.
		Then $\widehat W(z_0)=\int_0^\infty\varphi(x)\,\mathrm dx/x>0$.
		Thus $D$ is holomorphic near $z_0$.
		Dividing this function by the two nonzero factors above shows that $1/\zeta_\lambda$
		is holomorphic in that neighborhood.
		Finally, $1-2^{-z}$ has no zeros in this half-plane, so
		$\zeta_{\mathbb Q(i)}(z)$ has no zeros there.
	\end{proof}
	\section{AI statement}
	The author provided the idea, GPT6.1Sol produced the initial proof,
	and the author then improved that proof.
	The author changed the order of the material, removed cumbersome
	passages and notation, and expanded some points that were unclear.
	The author checked every paragraph and takes responsibility for the
	entire paper.

\end{document}